\documentstyle[12pt]{article}

\def\1{\'{\i}}

\def\be{\begin{equation}}
\def\ee{\end{equation}}
\def\bea{\begin{eqnarray}}
\def\eea{\end{eqnarray}}

\def\p{\varepsilon}

\def\jj{K} 
\def\pu{{H}} 
\def\pd{{P}}
\def\cu{{C_1}}
\def\cd{{C_2}}
\def\dd{D}

\def\te{\chi}

\def\sss{{\cal S}} 
\def\ggg{{\cal G}} 
 
\def\R{{\rm I\kern-.2em R}} 
 
\def\>#1{{\bf #1}}                 
\def\1{\'{\i}}                           
\def\R{\rm I\kern-.2em R}

\begin{document}

\ \vspace{2cm}

\begin{center}
{\LARGE{\bf{Null-plane Quantum Poincar\'e Algebras}}}

{\LARGE{\bf{and their Universal $R$-matrices}}}
\end{center}

\bigskip 

\begin{center}
A. Ballesteros${\dagger}$, F. J. Herranz${\dagger}$, M.A. del Olmo${\ddagger}$
and Mariano Santander${\ddagger}$
\end{center}

\begin{center}
\em 
{ { {}${\dagger}$ Departamento de F\1sica, Universidad de Burgos}
\\  E-09001, Burgos, Spain}

{ { {}${\ddagger}$ Departamento de F\1sica Te\'orica, Universidad de
Valladolid } \\   E-47011, Valladolid, Spain }
\end{center}
\rm

\begin{abstract}
A non-standard quantum deformation  of the  Poincar\'e algebra  is
presented in a null-plane framework for 1+1, 2+1 and 3+1 dimensions.
Their corresponding universal $R$-matrices are obtained in a
factorized form by choosing  suitable bases related to the $T$-matrix
formalism.
\end{abstract}

\section{Introduction}

 Among quantum deformations of the Poincar\'e algebra
we find three remarkable Hopf structures:

\noindent
$\bullet$ Two standard  deformations as quantizations of 
  coboundary Lie bialgebras coming from skew-solutions
of the modified classical Yang--Baxter equation (YBE). Physically these
  are obtained in a purely kinematical framework encoded within the usual
Poincar\'e basis. They are the
$\kappa$-Poincar\'e algebra
\cite{Luka,Gill,Lukb}  where the deformation parameter can  be
interpreted as a fundamental time scale  and a
$q$-Poincar\'e algebra \cite{Afin} where the quantum parameter is a
fundamental length.

\noindent
$\bullet$ One non-standard deformation as a quantization of
a (triangular) coboundary Lie  bialgebra coming from a skew-solution 
of the classical YBE.  It is closely related to the  
Jordanian or $h$-deformation for
$sl(2,\R)$ introduced in \cite{Demidov,Zakr,Ohn}. This structure 
  has been called    `null-plane'  quantum Poincar\'e
algebra   \cite{BeyondNull}  since it is constructed in a
null-plane context where the  Poincar\'e invariance splits into a kinematical
and dynamical part
\cite{Dirac,LS}.

  On the problem of obtaining universal $R$-matrices for the above
quantum Poincar\'e algebras  the results, as far as we know, are as
follows:

\noindent
$\bullet$ For  standard  Poincar\'e deformations,    a universal 
$R$-matrix has been only found    in  2+1 dimensions \cite{Itab} by means
of a contraction procedure starting from $so(4)_q$.  There are not 
known universal $R$-matrices  neither for the 1+1 case
 nor for the 3+1 case,  although for 1+1  there
are partial results \cite{soby}.

\noindent
$\bullet$ For  non-standard Poincar\'e  deformations  some results have been
obtained recently by   following
different methods:

\noindent
{\em 1+1 case}.  i) By an explicit  construction \cite{Irana,RR}
leading to a expression    similar    to the one obtained in
\cite{Vladi} for a Hopf subalgebra of $U_hsl(2,\R)$. ii) Following a universal 
$T$-matrix construction \cite{fg} giving rise to a factorized form
\cite{TT}. iii)   A contraction process from
the universal $R$-matrix of $U_hsl(2,\R)$ can also be applied.

\noindent
{\em 2+1  case}. By means of a contraction  method, taking as
starting point the non-standard deformation of $so(2,2)$
\cite{Iranb,poin}.

\noindent
{\em 3+1  case}. By applying a non-linear change of basis  
inspired in the results of the above cases and by using a universal
$T$-matrix approach; the resulting $R$-matrix  is 
factorized \cite{rnull}.

  Factorized  universal $R$-matrices    are given by
ordered (usual) exponentials of the elements appearing within the
corresponding classical $r$-matrix, thus they adopt a simple form.
This kind of factorized expressions are relevant
when $T$ is interpreted as a transfer matrix  in quantum field theory
\cite{Fr}. These $R$-matrices could be useful in order
to construct new integrable examples linked to the null-plane
evolution scheme (e.g., the infinite momentum frame
approach \cite{KSRep}, gauge field theory quantized on
null-planes \cite{KSop} and applications in Hadron
spectroscopy \cite{HHLS}).

In the sequel we present the results  concerning the factorized
universal $R$-matrices for  the non-standard quantum  Poincar\'e
algebra in a null-plane framework up to  3+1 
dimensions, paying special attention to the $T$-matrix method.  There is a 
main point in our procedure:    the Poincar\'e generators involved
in the non-standard classical $r$-matrix provide    after quantization a Hopf
subalgebra
$U_z\sss$    and it is possible to obtain its
corresponding  universal $R$-matrix, this is, a solution of the
quantum YBE which   fulfills the property
\be
{\cal R}\Delta (X) {\cal R}^{-1}=\sigma\circ \Delta (X),\qquad
 \forall X\in  U_z\sss.
\label{za}
\ee
The remarkable fact is that (\ref{za}) is also verified by the
Poincar\'e generators out of the Hopf subalgebra, thus ${\cal R}$ is
a universal $R$-matrix for the complete quantum Poincar\'e algebra.


\section {${U_z{\cal P}(1+1)}$: via universal $T$-matrix}

We recall first the main features of the $T$-matrix approach.
Let $U_z\ggg$ be a quantum algebra and $Fun_z(G)$ its associated dual
Hopf algebra or quantum group. The universal $T$-matrix of $U_z\ggg$
is the Hopf algebra dual form
$T=\sum_{\mu} X^{\mu}
\otimes p_{\mu}$ where $\{X^{\mu}\}$ is a basis for $U_z\ggg$
and $\{p_{\mu}\}$ its dual in $Fun_z(G)$, this is,
$\langle p_{\nu},X^{\mu}\rangle =\delta_{\nu}^{\mu}$  \cite{fg}. 
In general, the $T$-matrix can be interpreted as the universal
$R$-matrix for the quantum double linked to $U_z\ggg$
\cite{fdal,Chang,FRT} (within this quantum double, $U_z\ggg$
is a Hopf subalgebra). As a consequence, if there exists an algebra
isomorphism  and coalgebra anti-isomorphism 
$\Phi$ between   $U_z\ggg$ and   $Fun_z(G)$, the $T$-matrix can be
completely written in terms of the generators of $U_z\ggg$, in such a
way that
\be
{\cal R}=(\mbox{id}\otimes \Phi){\cal T},
\label{ab}
\ee
 with $\Phi$ acting on the generators of $Fun_z(G)$, is a universal
$R$-matrix for $U_z\ggg$.

We   illustrate now the above ideas
working out the `toy example' corresponding to the 1+1 Poincar\'e
algebra ${\cal P}(1+1)$. We choose as Poincar\'e generators the 
boost $K$, and the translations along the light-cone $P_+$,
$P_-$; they satisfy
\be
[K,P_+]=P_+,\qquad [K,P_-]=-P_-,\qquad [P_-,P_+]=0.
\label{ac}
\ee
The non-standard classical $r$-matrix   $r=2z\,K\wedge P_+$ 
provides  the  cocommutators  by means of
$\delta(X)=[1\otimes X + X \otimes 1,\,  r]$:
\be
\delta(P_+)=0,\quad \delta(K)=2z\, K\wedge P_+,
\quad \delta(P_-)=2z\, P_-\wedge P_+ .
\label{ad}
\ee
The coproduct and commutators of the Hopf algebra $U_z{\cal
P}(1+1)$ which deform  this Poincar\'e bialgebra are
\bea 
 &&\Delta (P_+)  =1 \otimes P_+  + P_+ \otimes 1,\quad
 \Delta (P_-) =1 \otimes P_- + P_- \otimes e^{2 z P_+ },\cr
&&  \Delta (K ) =1 \otimes K  + K  \otimes e^{2 z P_+ },
\label{ae}
\eea 
 \be
 [K ,P_+]=\frac 1{2z}(e^{2zP_+}  -1 ),\qquad 
[K ,P_-]=-  P_-  ,\qquad [P_+,P_-]=0.
\label{af}
\ee 

   Let us focus   on the Hopf subalgebra $U_z\sss$ generated by $K$
and $P_+$; its  quantum dual group $Fun_z(S)$ has coordinates $\hat
\te$ and $\hat a_+$ and its structure is given by:   
\be
\Delta({\hat \te})={\hat \te}\otimes 1 + 1 \otimes {\hat
\te},\quad \Delta({\hat a_+})= {\hat a_+}\otimes 1 + e^{\hat
\te}\otimes {\hat a_+}  ,\quad
[{{\hat \te}},{{\hat a_+}}]=2z  (e^{ {\hat \te}} -1).
\label{ag}
\ee 
The  associated universal $T$-matrix is \cite{TT}:
\be
{\cal T}= \exp\{ {P_+}\otimes {\hat a_+} \}\exp\{K
\otimes {\hat \te}\}.
\label{ah}
\ee
By taking into account (\ref{ae}), (\ref{af}) and (\ref{ag}) it can be
easily checked that the map  $\Phi$ defined by
\be
\Phi(\hat\te)=2zP_+,\qquad  \Phi(\hat a_+)=-2z K  ,\qquad \Phi(1)=1,
\label{ai}  
\ee
is an algebra isomorphism  and coalgebra anti-isomorphism 
  between   $U_z\sss$ and   $Fun_z(S)$, so that
\be
{\cal R}= ( id\otimes\Phi)\, T  =\exp\{-2zP_+\otimes K\}
\exp\{2zK\otimes P_+\},
\label{aj}
\ee
  is a universal $R$-matrix for $U_z\sss$. Moreover, it
can be proven that $P_-$ verifies (\ref{za}) with respect to
(\ref{aj}), hence we conclude that  this is a universal
$R$-matrix for the whole $U_z{\cal P}(1+1)$.


\section {$U_w{\cal P}(2+1)$: via contraction}

Instead of    developping    the $T$-matrix procedure,  we
show in this case another       method: the contraction process.
We take as the starting point  the coproduct, 
commutation relations and universal $R$-matrix of the non-standard
quantum  $U_z sl(2,\R)=\langle A,A_+,A_-\rangle$ 
\cite{rrr}:
\bea 
&&\Delta (A_+ ) =1 \otimes A_+  + A_+ \otimes 1,\quad
  \Delta (A) =1 \otimes A + A\otimes
e^{2z A_+ }, \cr
 &&\Delta (A_-) = 1 \otimes A_- + A_-\otimes
e^{2 z A_+ }, 
 \eea 
\be 
  [A,A_+ ]= \frac{e^{2 z A_+} -1  } z,\quad 
[A,A_-]=-2 A_- +z A^2,\quad  [A_+  ,A_- ]= A .
\label{aab} 
\ee 
\be
 R_z=\exp\{-z A_+\otimes A\}\exp\{z A\otimes A_+\}.
\label{aac}
\ee
At this point, notice that  $U_z{\cal P}(1+1)$ (and its universal
$R$-matrix) can be recovered by
contracting $U_z sl(2,\R)$ as the limit $\varepsilon\to 0$ of the 
mapping defined by: 
\be
K=A/2,\quad P_+=\varepsilon A_+,\quad P_-=\varepsilon A_-,\quad w=
z/\varepsilon .
\ee

A non-standard quantum Poincar\'e algebra $U_w{\cal P}(2+1)$ can be
obtained in two steps \cite{poin}:

\noindent
$\bullet$ Take two    quantum
$sl(2,\R)$ algebras:
$U_z^{(1)} sl(2,\R)=\langle A^1,A_+^1,A_-^1\rangle$ and 
$U_{-z}^{(2)} sl(2,\R)=\langle A^2,A_+^2,A_-^2\rangle$. In this way
the following  generators  
\be
\begin{array}{lll}
  \jj=\frac 12(A^1- A^2),&\quad \dd=\frac 12(A^1+A^2),&\quad
\cu=   - A_-^1 - A_-^2 ,\cr
  \pu=   A_+^1 + A_+^2 ,&\quad 
\pd= A_+^1 - A_+^2,&\quad \cd= A_-^1 -
A_-^2,
\end{array}
\label{aae}
\ee
give rise to $U_z so(2,2)\simeq  U_z^{(1)}  sl(2,{\R}) \oplus
U_{-z}^{(2)}  sl(2,\R)$ in a conformal basis.
 
\noindent
$\bullet$ Apply the contraction defined by:
\be
\begin{array}{lll}
{P_+}=\varepsilon \, \frac 1{\sqrt{2}} \pd,&\quad  
{P_1}=\varepsilon\, \jj,&\quad 
{P_-}=- \varepsilon\, \frac 1{\sqrt{2}} \cd,\cr
 {E_1}=-\frac 1{\sqrt{2}} \pu,&\quad 
{F_1}= \frac 1{\sqrt{2}} \cu,&\quad {K_2}=D,\qquad 
w= \frac 1{\varepsilon\sqrt{2}} z.
\end{array}
\label{ccb}
\ee

After the limit $\varepsilon\to 0$, the contracted generators
$\{P_+,P_-,P_1,E_1,F_1,K_2\}$ and deformation  parameter $w$ close a
2+1 quantum Poincar\'e algebra $U_w{\cal P}(2+1)$ in a null-plane
basis:
\be
\begin{array}{l}
  \Delta (P_+)  =1 \otimes P_+  + P_+ \otimes 1,\qquad
\Delta (E_1)  =1 \otimes E_1  + E_1 \otimes 1,\cr
  \Delta (P_-) =1 \otimes P_- + P_- \otimes e^{2 w P_+ },
\qquad \Delta (P_1) =1 \otimes P_1
+ P_1 \otimes e^{2 w P_+ },\cr
  \Delta (F_1) =1 \otimes F_1
+ F_1 \otimes e^{2 w P_+ } - 2 w P_-\otimes e^{2 w P_+ }E_1,
\cr
 \Delta (K_2) =1 \otimes K_2
+ K_2 \otimes e^{2 w P_+ } - 2 w P_1\otimes e^{2 w P_+ }E_1,
\end{array}
\label{ccc}
\ee
\be
\begin{array}{l}
 [K_2,P_+]=\frac 1{2w}(e^{2wP_+}  -1 ),\qquad 
[K_2,P_-]=- P_- - wP_1^2 ,\cr
 [K_2,E_1]=E_1 e^{2wP_+}  ,\qquad 
[K_2,F_1]=- F_1 - 2wP_1 K_2 , \cr
 [E_1,P_1]=\frac 1{2w}(e^{2wP_+}  -1 ),\qquad 
[F_1,P_1]=  P_- + wP_1^2 ,\cr
 [E_1,F_1]=K_2,\qquad 
[P_+,F_1]=  -  P_1,\qquad 
[P_-,E_1]=  -  P_1  .
\end{array}
\label{ccd}
\ee

The universal $R$-matrix for $U_w{\cal P}(2+1)$ also comes from the
contraction (\ref{ccb}) applied onto the one corresponding to $U_z
so(2,2)$ that is 
$R_z^{(1)}R_{-z}^{(2)}$:
\be 
{\cal R}_w=\exp\{2 w E_1 \otimes P_1  \}
\exp\{-2w   P_+ \otimes K_2 \}
\exp\{2w   K_2\otimes P_+  \}
\exp\{-2w  P_1\otimes E_1  \}.
\ee
The first order in $w$ gives the classical $r$-matrix  
$r=2w(K_2\wedge P_+ + E_1\wedge P_1)$ which provides the Lie
bialgebra underlying this quantum  Poincar\'e algebra.

  Obviously this method can be only carried out when the adequate
quantum semisimple algebra structure is known. This is not the case
  for next dimension so  we use again the $T$-matrix formalism.


\section {$U_z{\cal P}(3+1)$: via universal $T$-matrix}

Let us introduce first the classical Poincar\'e (bi)algebra in a
      null-plane basis.
We   consider the null-plane `orthogonal'
to the light-like vector $n=(\frac 12,0,0,\frac 12)$ as
our initial surface \cite{LS}. A coordinate  system
well adapted to null-planes  $\Pi_n^\tau$ is given by 
\be
x^- = n\cdot x =\frac 12(x^0-x^3)=\tau,\quad
x^+=x^0+x^3,\quad x_T=(x^1,x^2). \label{aa} 
\ee 
A point $x$ contained in the null-plane $\Pi_n^\tau$
is labelled by the coordinates $(x^+, x^1,x^2)$; the
remaining coordinate $x^-$ plays the role of an    evolution parameter
(`time').   
 A particular null-plane is $\Pi_n^0, (\tau=0)$, i.e., $x^- =
n\cdot x =0$ (it  is    invariant under the action of   the
boosts generated by $K_3$). A basis    of ${\cal P}(3+1)$  adapted to 
these coordinates  is $\{P_+,P_i,P_-,E_i,F_i,K_3,J_3\}$, where the   generators
$P_+$, $P_-$, $E_i$ and
$F_i$    are defined in the terms of the usual  kinematical
ones $\{P_0,P_j,K_j,J_j\}$ $(j=1,2,3)$ by:
\be
\begin{array}{lll}
{P_+}=\frac 12(P_0+P_3),&\qquad {P_-}=P_0-P_3,&\qquad 
{E_1}=\frac 12(K_1+J_2),\cr
 {F_1}=K_1-J_2,&\qquad {F_2}=K_2+J_1,&\qquad {E_2}=\frac 12(K_2-J_1).
\end{array}
\label{bbb}
\ee
  The non-vanishing Lie brackets of ${\cal P}(3+1)$ are
$(i,j=1,2)$: 
\be
\begin{array}{l}
 [K_3,P_+]=P_+,\quad
[K_3,P_-]=-P_-,\quad [K_3,E_i]=E_i,\quad 
[K_3,F_i]=-F_i,\cr
 [J_3,P_i]=-\p_{ij3}P_j, \qquad  [J_3,E_i]=-\p_{ij3}E_j,
\qquad [J_3,F_i]=-\p_{ij3}F_j,\cr 
[E_i,P_j]=\delta_{ij}P_+,\qquad  [F_i,P_j]=\delta_{ij}P_-,
\qquad [P_+,F_i]=-P_i,\cr 
[E_i,F_j]=\delta_{ij}K_3 +\p_{ij3}J_3,\qquad
[P_-,E_i]=-P_i.
\end{array}
\label{com2}
\ee
The stability group of the
plane $\Pi_n^0$ is generated by $\{P_+,P_i,E_i,K_3,J_3\}$. The
remaining three generators  act
on $\Pi_n^0$ as follows:   $P_-$ translates
$\Pi_n^0 $ into $\Pi_n^\tau$, while   both   $F_i$
rotate it around the surface of the light-cone $x^2=0$.
Therefore, if $x^-=\tau$ is considered as an evolution
parameter, then $P_-$ and $F_i$ describe the dynamical
evolution  from the null-plane $x^-=0$.

The classical $r$-matrix we choose is the natural generalization of
the one corresponding to the 2+1 case:
\be
r= 2z (K_3 \wedge P_+ +E_1\wedge P_1
+ E_2\wedge P_2) .
\label{bg}
\ee
Therefore the Poincar\'e cocommutators are:
\be
\begin{array}{l}
\delta(X)=0, \qquad \mbox{for}\quad
X\in\{P_+,E_i,J_3\},\cr 
 \delta(Y)=2z(Y\wedge P_+) ,\quad  
\mbox{for}\quad Y\in\{P_-,P_i\},\cr 
  \delta(F_1)=2z(F_1\wedge P_+ - P_-\wedge E_1 - P_2\wedge J_3),\cr
  \delta(F_2)=2z(F_2\wedge P_+ - P_-\wedge E_2 + P_1\wedge J_3),\cr
  \delta(K_3)=2z(K_3\wedge P_+ - P_1\wedge E_1 - P_2\wedge E_2).
\end{array}
\label{bbgg}
\ee
The quantum deformation $U_z{\cal P}(3+1)$of the Lie bialgebra so
obtained is given by \cite{rnull}:
\be
\begin{array}{l}
\Delta(X)=1\otimes X+X\otimes 1, \qquad \mbox{for}\quad
X\in\{P_+,E_i,J_3\},\cr 
 \Delta(Y)=1\otimes Y+Y\otimes e^{2zP_+},\quad  
\mbox{for}\quad Y\in\{P_-,P_i\},\cr 
  \Delta(F_1)=1\otimes F_1
+F_1\otimes e^{2zP_+} - 2zP_-\otimes E_1 e^{2zP_+}
 - 2z P_2\otimes J_3 e^{2zP_+},\cr
  \Delta(F_2)=1\otimes F_2+F_2\otimes e^{2zP_+}
 - 2zP_-\otimes E_2 e^{2zP_+} + 
2z P_1\otimes J_3 e^{2zP_+},\cr 
 \Delta(K_3)=1\otimes K_3+K_3\otimes e^{2zP_+}
 - 2zP_1\otimes E_1 e^{2zP_+} - 
2z P_2\otimes E_2 e^{2zP_+};
\end{array}
\label{bb}
\ee
\be
\begin{array}{l}
  [K_3,P_+]=\frac 1{2z} ({e^{2zP_+} - 1}),\qquad [K_3,P_-]=-P_-
-zP_1^2- zP_2^2,\cr 
  [K_3,E_i]=E_ie^{2zP_+},\qquad [K_3,F_i]=- F_i - 2 z K_3 P_i,\cr 
  [J_3,P_i]=-\p_{ij3}P_j, \qquad [J_3,E_i]=-\p_{ij3}E_j,
\qquad [J_3,F_i]=-\p_{ij3}F_j,\cr
  [E_i,P_j]=\delta_{ij}\frac 1{2z} (  {e^{2zP_+} - 1}),
\qquad  [F_i,P_j]=\delta_{ij}(P_- + zP_1^2+ zP_2^2),\cr
  [E_i,F_j]=\delta_{ij}K_3 +\p_{ij3}J_3e^{2zP_+},\qquad
[P_+,F_i]=-P_i,\cr 
  [F_1,F_2]=2z(P_1F_2 - P_2F_1), \qquad [P_-,E_i]=-P_i .
\end{array}
\label{bc}
\ee

The six generators appearing in the classical $r$-matrix (\ref{bg})
  close a quantum Hopf subalgebra $U_z \sss$.
The universal $T$-matrix for this Hopf subalgebra can be computed
and reads \cite{rnull}:
\be
{\cal T}=e^{E_2\otimes{\hat e}_2}e^{E_1\otimes{\hat e}_1}
e^{P_+\otimes{\hat a}_+}e^{K_3\otimes{\hat k}_3}
e^{P_1\otimes{\hat a}_1}e^{P_2\otimes{\hat a}_2},
\label{a14}
\ee
where ${\hat a}_+,{\hat a}_i,{\hat e}_i,{\hat k}_3$ are dual to the 
corresponding Poincar\'e generators. This canonical element leads to
a  factorized  universal $R$-matrix for $U_z \sss$:
\bea
&&{\cal R}=\exp\{2 z E_2\otimes P_2\}\exp\{2 z E_1\otimes P_1\}
\exp\{- 2 z P_+\otimes K_3\}\cr
&&\qquad \times \exp\{2 z K_3\otimes P_+\}
\exp\{-2 z P_1\otimes E_1\}\exp\{- 2 z P_2\otimes  E_2\} .
\label{cq}
\eea 
 Furthermore it can be shown  that ${\cal R}$
satisfies the property (\ref{za}) for  the four remaining generators:
$P_-$, $F_1$, $F_2$ and $J_3$. Therefore, it is  a
universal $R$-matrix for the whole null-plane quantum Poincar\'e
algebra.


\section*{Acknowledgments}
This work has been
partially supported by the Spanish DGICYT (Projects PB94--1115).



\begin{thebibliography}{99}
\itemsep=-.2pc


\bibitem{Luka}
J.~Lukierski, A.~Nowicky,  H.~Ruegg and V.~Tolstoy,      
  {\em Phys. Lett. B}   {\bf  264}  (1991)  331.


\bibitem{Gill}
S.~Giller, P.~Koshinski, J.~Kunz, M.~Majewski and P.~Maslanka,
  {\em Phys. Lett. B}   {\bf  286}  (1992)   57.


\bibitem{Lukb}
J.~Lukierski,   A.~Nowicky and H.~Ruegg,      
  {\em Phys. Lett. B}   {\bf  293}  (1992)  344.

 
 
\bibitem{Afin}
A.~Ballesteros,  F.J.~Herranz, M.A.~del Olmo and 
M.~Santander, {\em J. Math. Phys.}  {\bf 35}  (1994)  4928.

 

\bibitem{Demidov}
E.E.~Demidov, Y.I.~Manin, E.E.~Mukhin and D.V.~Zhdanovich,
{\em Progr. Theor. Phys. Suppl.} {\bf 102} (1990) 203.

\bibitem{Zakr}
S.~Zakrzewski, {\em Lett. Math. Phys.} {\bf 22} (1991) 287.

 
\bibitem{Ohn} C.~Ohn,
{\em Lett. Math. Phys.} {\bf 25} (1992) 85.

 
\bibitem{BeyondNull}
A.~Ballesteros,  F.J.~Herranz, M.A.~del Olmo  and M.~Santander,  
{\em  J. Phys. A: Math. Gen.}  {\bf 28}  (1995) 941; {\em  Phys. Lett. B}  
{\bf  351}  (1995) 137.


\bibitem{Dirac} P.A.M.~Dirac, {\em Rev. Mod. Phys.} 
{\bf 21} (1949) 392.


\bibitem{LS} H.~Leutwyler and J.~Stern, 
{\em Ann. Phys. (N.Y.)}  {\bf 112}   (1978)  94.


\bibitem{Itab}
E.~Celeghini, R.~Giachetti, E.~Sorace and M.~Tarlini, 
 ``Contractions of quantum groups",  in:
{\em Lecture Notes in Math.~n.~1510}, p.~221,
(Springer-Verlag, Berl\1n 1992).
 
\bibitem{soby}
J.~Sobczyk,   
{\em  J. Phys. A: Math. Gen.}  {\bf 29}  (1996) 2887.


\bibitem{Irana}
M.~Khorrami, A.~Shariati, M.R.~Abolhassani and 
A.~Aghamohammadi, 
 {\em Mod. Phys. Lett. A}  {\bf 10}  (1995)  873.


\bibitem{RR}
A.~Ballesteros, E.~Celeghini, F.J.~Herranz, M.A.~del 
Olmo and M.~Santander,  {\em J. Phys. A: Math. Gen.}  {\bf 28} 
(1995)  3129.

\bibitem{Vladi}
A.A.~Vladimirov, 
 {\em Mod. Phys. Lett. A}  {\bf 8}  (1993)  2573.



\bibitem{fg} C.~Fronsdal and A.~Galindo, 
{\em Lett. Math. Phys.}
{\bf 27}  (1993) 39.

 
\bibitem{TT}
A.~Ballesteros,  F.J.~Herranz, M.A.~del~Olmo, 
C.M.~Pere\~na and M.~Santander,   {\em J. Phys. A: Math. Gen.}  {\bf
28}  (1995) 7113.

 
\bibitem{Iranb}
A.~Shariati, A.~Aghamohammadi  and  M.~Khorrami, 
 {\em Mod. Phys. Lett. A}  {\bf  11} (1996)  187.

 \bibitem{poin}
A.~Ballesteros  and  F.J.~Herranz, 
 {\em Mod. Phys. Lett. A}  {\bf 11} (1996) 1745.


 
 \bibitem{rnull}
A.~Ballesteros,  F.J.~Herranz and C.M.~Pere\~na, ``Null-plane Quantum
Universal $R$-matrix", {\em Phys. Lett. B}  in press.


 
\bibitem{Fr} C.~Fronsdal, {\em Lett. Math. Phys.}
{\bf 24}  (1992) 73.

 
\bibitem{KSRep} I.B.~Kogut and L.~Susskind, 
{\em  Phys. Reports C}  {\bf 8}  (1973) 75.

\bibitem{KSop} I.B.~Kogut and D.E.~Soper, 
{\em  Phys. Rev. D} {\bf 1} (1970) 2901.


\bibitem{HHLS} H.H.~Liu and D.E.~Soper, {\em  Phys. Rev. D}
{\bf 48}  (1993) 1841.


\bibitem{fdal}  C.~Fronsdal,  ``Universal $T$-Matrix for
Twisted Quantum $gl(N)$", preprint UCLA/93/TEP/3 (1993). 


\bibitem{Chang} Z.~Chang,  {\em Phys. Rep.} {\bf 262} (1995)  137. 


\bibitem {FRT} 
 N.Y.~Reshetikhin, L.A.~Takhtadzhyan  and  L.D.~Faddeev,   
{\em Leningrad Math. J.} {\bf 1}   (1990) 193.   


 \bibitem{rrr}
A.~Ballesteros  and  F.J.~Herranz, 
{\em J. Phys. A: Math. Gen.}  {\bf
29}  (1996) L311.








 
\end{thebibliography}
\end{document}